\documentclass[11pt, twoside]{article}
\usepackage{latexsym}
\usepackage{amsmath}
\usepackage{amssymb}
\usepackage[all]{xy}
\usepackage{amsfonts}
\usepackage{verbatim}
\usepackage{amsthm}
\usepackage{bm}
\usepackage{mathrsfs}
\usepackage{epsfig}
\usepackage{xy}
\usepackage{array}
\usepackage{stmaryrd}
\usepackage{color}
\usepackage[colorlinks=true,linkcolor=blue,citecolor=blue]{hyperref}
\usepackage{tikz}
\usetikzlibrary{arrows,calc}
\usepackage{etex}
\usepackage{mathdots}
\usepackage{float}
\usepackage{graphics}
\usepackage{pdflscape}
\usepackage{CJK}
\usepackage{anysize,hyperref}
\input xypic
\xyoption{all}
\usepackage{extarrows}
\usepackage[perpage,symbol]{footmisc}
\usepackage{setspace}
\def\E{\mathbb{E}}

\def\C{\mathcal{C}}
\def\s{\mathfrak{s}}

\def\del{\delta}
\def\dr{\ar@{->}[r]}

\begin{document}
\baselineskip=15pt
\title{\Large{\bf Homotopy cartesian squares in extriangulated categories$^\bigstar$\footnotetext{\hspace{-1em}$^\bigstar$Panyue Zhou was supported by the National Natural Science Foundation of China (Grant No. 11901190).}}}
\medskip
\author{Jing He, Chenbei Xie and Panyue Zhou}

\date{}

\maketitle
\def\blue{\color{blue}}
\def\red{\color{red}}

\newtheorem{theorem}{Theorem}[section]
\newtheorem{lemma}[theorem]{Lemma}
\newtheorem{corollary}[theorem]{Corollary}
\newtheorem{proposition}[theorem]{Proposition}
\newtheorem{conjecture}{Conjecture}
\theoremstyle{definition}
\newtheorem{definition}[theorem]{Definition}
\newtheorem{question}[theorem]{Question}
\newtheorem{remark}[theorem]{Remark}
\newtheorem{remark*}[]{Remark}
\newtheorem{example}[theorem]{Example}
\newtheorem{example*}[]{Example}
\newtheorem{condition}[theorem]{Condition}
\newtheorem{condition*}[]{Condition}
\newtheorem{construction}[theorem]{Construction}
\newtheorem{construction*}[]{Construction}

\newtheorem{assumption}[theorem]{Assumption}
\newtheorem{assumption*}[]{Assumption}

\baselineskip=17pt
\parindent=0.5cm

\begin{abstract}
\baselineskip=16pt
Let $(\C,\E,\s)$ be an extriangulated category. Given a composition of two commutative squares in $\C$, if two commutative squares are homotopy cartesian, then their composition is also a homotopy cartesian.
This covers the result by Mac Lane (1998) for abelian categories and the result by Christensen and Frankland (2022) for triangulated categories.\\[0.5cm]
\medskip
\end{abstract}

\pagestyle{myheadings}
\markboth{\rightline {\scriptsize   J. He, C. Xie and P. Zhou}}
         {\leftline{\scriptsize  Homotopy cartesian squares in extriangulated categories}}

\section{Introduction}
Most of homological algebra can be carried out in the setting of abelian
categories, a class of categories which includes all categories of modules.
We recall the properties of pull-back and push-out squares in an abelian category.
Given two morphisms $x'\colon A'\to B'$  and $b\colon B\to B'$. Then a commutative diagram
$$\xymatrix{A\ar[d]_{a}\ar[r]^{x} &B\ar[d]^{b}\\
 A'\ar[r]^{x'} & B'}$$
is \emph{pull-back} if and only if
$0\xrightarrow{~}A\xrightarrow{~\left[
              \begin{smallmatrix}
                a\\ x
              \end{smallmatrix}
            \right]~}A'\oplus B\xrightarrow{~\left[
              \begin{smallmatrix}
            x'&-b
              \end{smallmatrix}
            \right]~}B'$
is left exact.  Given two morphisms $x\colon A\to B$ and $a\colon A\to A'$. Then a commutative diagram
$$\xymatrix{A\ar[d]_{a}\ar[r]^{x} &B\ar[d]^{b}\\
 A'\ar[r]^{x'} & B'}$$
is \emph{push-out} if and only if
$A\xrightarrow{~\left[
              \begin{smallmatrix}
                a\\ x
              \end{smallmatrix}
            \right]~}A'\oplus B\xrightarrow{~\left[
              \begin{smallmatrix}
            x'&-b
              \end{smallmatrix}
            \right]~}B'\xrightarrow{~}0$
is right exact. By the above properties, we can obtain that the commutative diagram
$$\xymatrix{A\ar[d]_{a}\ar[r]^{x} &B\ar[d]^{b}\\
 A'\ar[r]^{x'} & B'}$$
is \emph{pull-back and push-out} if and only if
$$0\xrightarrow{~}A\xrightarrow{~\left[
              \begin{smallmatrix}
                a\\ x
              \end{smallmatrix}
            \right]~}A'\oplus B\xrightarrow{~\left[
              \begin{smallmatrix}
            x'&-b
              \end{smallmatrix}
            \right]~}B'\xrightarrow{~}0$$
is exact.

Triangulated categories were introduced in the mid 1960's by Verdier
in his thesis. Having their origins in algebraic geometry and
algebraic topology, triangulated categories have by now become indispensable
in many different areas of mathematics.
Neeman \cite[Definition 1.4.1]{N} introduced the notion of homotopy cartesian squares in triangulated categories.
More precisely, let $\C$ be a triangulated category with shift functor $\Sigma$. A commutative square
$$\xymatrix{A\ar[d]_{a}\ar[r]^{x} &B\ar[d]^{b}\\
 A'\ar[r]^{x'} & B'}$$
s called \emph{homotopy cartesian} if there exists a triangle
$$A\xrightarrow{~\left[
              \begin{smallmatrix}
                a\\ x
              \end{smallmatrix}
            \right]~}A'\oplus B\xrightarrow{~\left[
              \begin{smallmatrix}
            x'&-b
              \end{smallmatrix}
            \right]~}B'\xrightarrow{~\partial~}\Sigma A$$
The morphism $\partial$ is called a \emph{differential }of the homotopy cartesian square.
We can think of a homotopy cartesian square as the triangulated analogue of
a pull-back and push-out square in an abelian category.

Recently, extriangulated categories were recently introduced by Nakaoka and Palu \cite{NP} by extracting those properties of ${\rm Ext}^1(-,-)$ on exact categories and on triangulated categories that seem relevant from the point of view of cotorsion pairs. In particular, triangulated categories and exact categories are extriangulated categories. There are a lot of examples of extriangulated categories which are neither triangulated categories nor exact categories, see \cite{NP,ZZ,HZZ,ZhZ,NP1}. Hence, many results on exact categories (abelian categories are exact categories) and triangulated categories can be unified in the same framework.
Motivated by this idea, we study
homotopy cartesian squares in an extriangulated category.
Let us recall the definition of homotopy cartesian squares in an extriangulated category.
\begin{definition}\label{y2}\cite[Definition 3.1]{H}
Let $(\C,\E,\s)$ be an extriangulated category.
Then a commutative square
$$\xymatrix{A\ar[d]_{a}\ar[r]^{x} &B\ar[d]^{b}\\
 A'\ar[r]^{x'} & B'}$$
is called \emph{homotopy cartesian} if there exists an $\E$-triangle
$$A\xrightarrow{~\left[
              \begin{smallmatrix}
                a\\ x
              \end{smallmatrix}
            \right]~}A'\oplus B\xrightarrow{~\left[
              \begin{smallmatrix}
            x'&-b
              \end{smallmatrix}
            \right]~}B'\overset{\theta}{\dashrightarrow}$$
The $\E$-extension $\theta\in\E(B',A)$ is called a \emph{differential} of the homotopy cartesian square.
\end{definition}
This gives a simultaneous generalization of push-out and pullback squares in an abelian category
and homotopy cartesian squares in triangulated categories.

The main result
of this article is the following, which is a simultaneous generalization of a result of
 Mac Lane \cite{M} for abelian categories and of a result of Christensen-Frankland \cite{CF} for triangulated categories.

\begin{theorem}\label{main1}{\rm (see Theorem \ref{main1} for more details)}
Let $(\C,\E,\s)$ be an extriangulated category.
Consider the following commutative diagram
$$\xymatrix{
A\ar[d]_{a} \ar[r]^{x} \ar@{}[dr]|{\rm (I)} & B \ar[d]^b\ar[r]^y\ar@{}[dr]|{\rm (II)}&C\ar[d]^c\\
A' \ar[r]_{x'} & B'\ar[r]_{y'}&C'
} $$
If  the two squares {\rm (I)} and {\rm (II)} are homotopy cartesian in $\C$,
 then the outside rectangle is a homotopy cartesian in $\C$.

Moreover, given differentials $\theta_L\in\E(B',A)$ and $\theta_R\in\E(C',B)$ of the left square and right square respectively, there exists a differential $\theta_P$ for the outside rectangle satisfying
$$\theta_R=x_{\ast}\theta_P~~\mbox{and}~~\theta_L=(y')^{\ast}\theta_P.$$
\end{theorem}
This article is organized as follows. In Section 2, we recall the definition
of an extriangulated category. In Section 3, we give a composition of two commutative squares, if both squares are homotopy cartesian, then their composition is also a homotopy cartesian.

\section{Preliminaries}
Let us briefly recall some definitions and basic properties of extriangulated categories from \cite{NP}.
We omit some details here, but the reader can find
them in \cite{NP}.

Let $\mathcal{C}$ be an additive category equipped with an additive bifunctor
$$\mathbb{E}: \mathcal{C}^{\rm op}\times \mathcal{C}\rightarrow {\rm Ab},$$
where ${\rm Ab}$ is the category of abelian groups. For any objects $A, C\in\mathcal{C}$, an element $\delta\in \mathbb{E}(C,A)$ is called an $\mathbb{E}$-extension.
Let $\mathfrak{s}$ be a correspondence which associates an equivalence class $$\mathfrak{s}(\delta)=\xymatrix@C=0.8cm{[A\ar[r]^x
 &B\ar[r]^y&C]}$$ to any $\mathbb{E}$-extension $\delta\in\mathbb{E}(C, A)$. This $\mathfrak{s}$ is called a {\it realization} of $\mathbb{E}$, if it makes the diagrams in \cite[Definition 2.9]{NP} commutative.
 A triplet $(\mathcal{C}, \mathbb{E}, \mathfrak{s})$ is called an {\it extriangulated category} if it satisfies the following conditions.
\begin{enumerate}
\item[\rm (1)] $\mathbb{E}\colon\mathcal{C}^{\rm op}\times \mathcal{C}\rightarrow \rm{Ab}$ is an additive bifunctor.

\item[\rm (2)] $\mathfrak{s}$ is an additive realization of $\mathbb{E}$.

\item[\rm (3)] $\mathbb{E}$ and $\mathfrak{s}$  satisfy the compatibility conditions in \cite[Definition 2.12]{NP}.
 \end{enumerate}

\begin{remark}
Note that both exact categories and triangulated categories are extriangulated categories, see \cite[Example 2.13]{NP} and extension closed subcategories of extriangulated categories are
again extriangulated, see \cite[Remark 2.18]{NP}. Moreover, there are extriangulated categories which
are neither exact categories nor triangulated categories, see \cite[Proposition 3.30]{NP}, \cite[Example 4.14]{ZZ}, \cite[Remark 3.3]{HZZ},\cite[Lemma 20]{ZhZ} and \cite[Theorem 4.22]{NP1}.
\end{remark}

We collect the following terminology from \cite{NP}.

\begin{definition}\label{y1}\cite[Definition 2.19]{NP}
Let $(\C,\E,\s)$ be an extriangulated category.
\begin{itemize}
\item[(1)] A sequence $A\xrightarrow{~x~}B\xrightarrow{~y~}C$ is called a {\it conflation} if it realizes some $\E$-extension $\del\in\E(C,A)$.
    In this case, $x$ is called an {\it inflation} and $y$ is called a {\it deflation}.

\item[(2)] If a conflation  $A\xrightarrow{~x~}B\xrightarrow{~y~}C$ realizes $\delta\in\mathbb{E}(C,A)$, we call the pair $( A\xrightarrow{~x~}B\xrightarrow{~y~}C,\delta)$ an {\it $\E$-triangle}, and write it in the following way.
$$A\overset{x}{\longrightarrow}B\overset{y}{\longrightarrow}C\overset{\delta}{\dashrightarrow}$$
We usually do not write this $``\delta"$ if it is not used in the argument.

\item[(3)] Let $A\overset{x}{\longrightarrow}B\overset{y}{\longrightarrow}C\overset{\delta}{\dashrightarrow}$ and $A^{\prime}\overset{x^{\prime}}{\longrightarrow}B^{\prime}\overset{y^{\prime}}{\longrightarrow}C^{\prime}\overset{\delta^{\prime}}{\dashrightarrow}$ be any pair of $\E$-triangles. If a triplet $(a,b,c)$ realizes $(a,c)\colon\delta\to\delta^{\prime}$, then we write it as
$$\xymatrix{
A \ar[r]^x \ar[d]^a & B\ar[r]^y \ar[d]^{b} & C\ar@{-->}[r]^{\del}\ar[d]^c&\\
A'\ar[r]^{x'} & B' \ar[r]^{y'} & C'\ar@{-->}[r]^{\del'} &}$$
and call $(a,b,c)$ a {\it morphism of $\E$-triangles}.
\end{itemize}
\end{definition}

\section{Proof of the Main Result}

First we prove the following crucial lemma.

\begin{lemma}\label{y3}
Let $(\C,\E,\s)$ be an extriangulated category.
Consider a homotopy cartesian square
$$\xymatrix{A\ar[d]^{a}\ar[r]^{x} & B\ar[d]^{b} &&\\
 A'\ar[r]^{x'} & B'\ar[r]^{y'} & C'\ar@{-->}[r]^{\delta'}&}$$
with given differential $\theta\in\E(B',A)$ and the bottom row is an $\E$-triangle in $\C$.
Then the square extends to a morphism of $\E$-triangles of the form
$(a,b,1)$, as illustrated in the diagram
$$\xymatrix{A\ar[d]^{a}\ar[r]^{x} & B\ar[r]^y\ar[d]^{b} &
C'\ar@{=}[d]\ar@{-->}[r]^{\delta}&\\
 A'\ar[r]^{x'} & B'\ar[r]^{y'} & C'\ar@{-->}[r]^{\delta'}&,}$$
satisfying $\theta=(y')^{\ast}\delta$.
\end{lemma}

\proof Since the square $$\xymatrix{A\ar[d]_{a}\ar[r]^{x} &B\ar[d]^{b}\\
 A'\ar[r]^{x'} & B'}$$
is homotopy cartesian, then there exists an $\E$-triangle
$$A\xrightarrow{~\left[
              \begin{smallmatrix}
                a\\ x
              \end{smallmatrix}
            \right]~}A'\oplus B\xrightarrow{~\left[
              \begin{smallmatrix}
            x'&-b
              \end{smallmatrix}
            \right]~}B'\overset{\theta}{\dashrightarrow}$$
By the dual of Proposition 3.17 in \cite{NP}, we obtain the following commutative diagram
of $\E$-triangles
$$\xymatrix@C=1cm{&A\ar[d]^{\left[
              \begin{smallmatrix}
          a\\x
              \end{smallmatrix}
            \right]}\ar@{=}[r]&A\ar[d]^x&\\
A'\ar@{=}[d]\ar[r]^{\left[
              \begin{smallmatrix}
            1\\0
              \end{smallmatrix}
            \right]\quad} &A'\oplus B\ar[r]^{\quad\left[
              \begin{smallmatrix}
          0&1
              \end{smallmatrix}
            \right]}\ar[d]^{\left[
              \begin{smallmatrix}
            x'&-b
              \end{smallmatrix}
            \right]} &
B\ar[d]^{-y}\ar@{-->}[r]^{0}&\\
 A'\ar[r]^{x'} & B'\ar@{-->}[d]^{\theta}\ar[r]^{y'} & C'\ar@{-->}[d]^{
  -\delta}\ar@{-->}[r]^{\delta'}&\\
 &&&}$$
which satisfies $\theta=(y')^{\ast}\delta$.
From the diagram, we have $y=y'b$. The commutative diagram
$$\xymatrix{A\ar@{=}[d]\ar[r]^{x} & B\ar[r]^{-y}\ar@{=}[d]&
C'\ar[d]^{-1}\ar@{-->}[r]^{-\delta}&\\
 A\ar[r]^{x} & B\ar[r]^{y} & C\ar@{-->}[r]^{\delta}&}$$
shows that
$\xymatrix{ A\ar[r]^{x} & B\ar[r]^{y} & C\ar@{-->}[r]^{\delta}&}$
is an $\E$-triangle.
Thus we obtain the following morphism of $\E$-triangles
$$\xymatrix{A\ar[d]^{a}\ar[r]^{x} & B\ar[r]^y\ar[d]^{b} &
C'\ar@{=}[d]\ar@{-->}[r]^{\delta}&\\
 A'\ar[r]^{x'} & B'\ar[r]^{y'} & C'\ar@{-->}[r]^{\delta'}&.}$$
This completes the proof.    \qed
\medskip

We now prove our main result.

\begin{theorem}\label{main1}
Let $(\C,\E,\s)$ be an extriangulated category.
Consider the following commutative diagram
$$\xymatrix{
A\ar[d]_{a} \ar[r]^{x} \ar@{}[dr]|{\rm (I)} & B \ar[d]^b\ar[r]^y\ar@{}[dr]|{\rm (II)}&C\ar[d]^c\\
A' \ar[r]_{x'} & B'\ar[r]_{y'}&C'
} $$
If  the two squares {\rm (I)} and {\rm (II)} are homotopy cartesian in $\C$,
 then the outside rectangle is a homotopy cartesian in $\C$.

Moreover, given differentials $\theta_L\in\E(B',A)$ and $\theta_R\in\E(C',B)$ of the left square and right square respectively, there exists a differential $\theta_P$ for the outside rectangle satisfying
$$\theta_R=x_{\ast}\theta_P~~\mbox{and}~~\theta_L=(y')^{\ast}\theta_P.$$
\end{theorem}

\proof Since the square {\rm (I)} is homotopy cartesian, then there exists an $\E$-triangle
$$A\xrightarrow{~\left[
              \begin{smallmatrix}
                a\\ x
              \end{smallmatrix}
            \right]~}A'\oplus B\xrightarrow{~\left[
              \begin{smallmatrix}
            x'&-b
              \end{smallmatrix}
            \right]~}B'\overset{\theta_L}{\dashrightarrow}$$
with differential $\theta_L\in\E(B',A)$.
Note that the following diagram
$$\xymatrix@C=1.2cm@R=1.3cm{A\ar@{=}[d]\ar[r]^{\left[
              \begin{smallmatrix}
          a\\x
              \end{smallmatrix}
            \right]\quad} &A'\oplus B\ar[r]^{\quad\left[
              \begin{smallmatrix}
            x'&-b
              \end{smallmatrix}
            \right]}\ar[d]^{\left[
              \begin{smallmatrix}
          0&1\\
          1&0
              \end{smallmatrix}
            \right]} &
B'\ar[d]^{-1}\ar@{-->}[r]^{\theta_L}&\\
 A\ar[r]^{\left[
              \begin{smallmatrix}
          x\\a
              \end{smallmatrix}
            \right]\quad} & B\oplus A'\ar[r]^{\quad\left[
              \begin{smallmatrix}
            b&-x'
              \end{smallmatrix}
            \right]} & B'\ar@{-->}[r]^{\theta_L}&}$$
is commutative and all vertical morphisms are isomorphisms,
thus we have that
\begin{equation}\label{t1}
\begin{array}{l}
A\xrightarrow{~\left[
              \begin{smallmatrix}
                x\\ a
              \end{smallmatrix}
            \right]~}B\oplus A'\xrightarrow{~\left[
              \begin{smallmatrix}
            b&-x'
              \end{smallmatrix}
            \right]~}B'\overset{\theta_L}{\dashrightarrow}
\end{array}
\end{equation}
is also an $\E$-triangle.
By the direct sum of the $\E$-triangle (\ref{t1}) and
 the trivial  $\E$-triangle
$0\overset{}{\longrightarrow}C\xrightarrow{~-1~}C\overset{0}{\dashrightarrow}$,
we obtain that
$$A\xrightarrow{~\left[
              \begin{smallmatrix}
                x\\ a\\0
              \end{smallmatrix}
            \right]~}B\oplus A'\oplus C\xrightarrow{~\left[
              \begin{smallmatrix}
            b&-x'&0\\
            0&0&-1
              \end{smallmatrix}
            \right]~}B'\oplus C\xymatrix@C=1.2cm{\ar@{-->}[r]^{\left[
              \begin{smallmatrix}
        \theta_L&0
              \end{smallmatrix}
            \right]\;\;}&}$$
is an $\E$-triangle.
Note that the following diagram
$$\xymatrix@C=1.8cm@R=1.6cm{A\ar@{=}[d]\ar[r]^{\left[
              \begin{smallmatrix}
          x\\a\\ 0
              \end{smallmatrix}
            \right]\qquad} &B\oplus A'\oplus C\ar[r]^{\qquad\left[
              \begin{smallmatrix}
            b&-x'&0\\
            0&0&-1
              \end{smallmatrix}
            \right]\quad}\ar[d]^{\left[
              \begin{smallmatrix}
          1&0&0\\
          0&1&0\\
          y&0&1
              \end{smallmatrix}
            \right]} &
B'\oplus C\ar@{=}[d]\ar@{-->}[r]^{\quad\left[
              \begin{smallmatrix}
        \theta_L&0
              \end{smallmatrix}
            \right]}&\\
 A\ar[r]^{\left[
              \begin{smallmatrix}
          x\\a\\yx
              \end{smallmatrix}
            \right]\qquad} & B\oplus A'\oplus C\ar[r]^{\quad\left[
              \begin{smallmatrix}
            b&-x'&0\\
            y&0&-1
              \end{smallmatrix}
            \right]} & B'\oplus C\ar@{-->}[r]^{\quad\left[
              \begin{smallmatrix}
        \theta_L&0
              \end{smallmatrix}
            \right]}&}$$
is commutative and all vertical morphisms are isomorphisms,
thus we have that
$$A\xrightarrow{~\left[
              \begin{smallmatrix}
                x\\ a\\yx
              \end{smallmatrix}
            \right]~}B\oplus A'\oplus C\xrightarrow{~\left[
              \begin{smallmatrix}
            b&-x'&0\\
            y&0&-1
              \end{smallmatrix}
            \right]~}B'\oplus C\xymatrix@C=1.2cm{\ar@{-->}[r]^{\left[
              \begin{smallmatrix}
        \theta_L&0
              \end{smallmatrix}
            \right]\;\;}&}$$
is also an $\E$-triangle.
This shows that the square
$$\xymatrix{A\ar[d]_x\ar[r]^{\left[
              \begin{smallmatrix}
              a\\yx
              \end{smallmatrix}
            \right]\quad} &A'\oplus C\ar[d]^{\left[
              \begin{smallmatrix}
                x'&0\\
                0&1
              \end{smallmatrix}
            \right]}\\
 B\ar[r]_{\left[
              \begin{smallmatrix}
                b\\y
              \end{smallmatrix}
            \right]\quad} & B'\oplus C}$$
is homotopy cartesian.

Since the square {\rm (II)} is homotopy cartesian, then there exists an $\E$-triangle
$$B\xrightarrow{~\left[
              \begin{smallmatrix}
                b\\ y
              \end{smallmatrix}
            \right]~}B'\oplus C\xrightarrow{~\left[
              \begin{smallmatrix}
            y'&-c
              \end{smallmatrix}
            \right]~}C'\overset{\theta_R}{\dashrightarrow}$$
with differential $\theta_R\in\E(C',B)$.
Consider the following commutative diagram
$$\xymatrix@C=1.6cm{A\ar[d]^x\ar[r]^{\left[
              \begin{smallmatrix}
          a\\ yx
              \end{smallmatrix}
            \right]\qquad} &A'\oplus C\ar[d]^{\left[
              \begin{smallmatrix}
          x'&0\\
          0&1
              \end{smallmatrix}
            \right]} &&\\
 B\ar[r]^{\left[
              \begin{smallmatrix}
        b\\y
              \end{smallmatrix}
            \right]\qquad} & B'\oplus  C\ar[r]^{\quad\left[
              \begin{smallmatrix}
            y'&-c
              \end{smallmatrix}
            \right]} & C'\ar@{-->}[r]^{\quad \theta_R}&,}$$
by Lemma \ref{y3}, we have the following morphism of $\E$-triangles
\begin{equation}\label{t2}
\begin{array}{l}
\xymatrix@C=1.6cm@R=1.4cm{A\ar[d]^x\ar[r]^{\left[
              \begin{smallmatrix}
          a\\ yx
              \end{smallmatrix}
            \right]\qquad} &A'\oplus C\ar[r]^{\qquad\left[
              \begin{smallmatrix}
         y'x'&-c
              \end{smallmatrix}
            \right]\quad}\ar[d]^{\left[
              \begin{smallmatrix}
          x'&0\\
          0&1
              \end{smallmatrix}
            \right]} &
C'\ar@{=}[d]\ar@{-->}[r]^{\quad \theta_P}&\\
 B\ar[r]^{\left[
              \begin{smallmatrix}
        b\\y
              \end{smallmatrix}
            \right]\qquad} & B'\oplus  C\ar[r]^{\quad\left[
              \begin{smallmatrix}
            y'&-c
              \end{smallmatrix}
            \right]} & C'\ar@{-->}[r]^{\quad \theta_R}&.}
\end{array}
\end{equation}
It follows that the $\E$-triangle
$$A\xrightarrow{~\left[
              \begin{smallmatrix}
                a\\ yx
              \end{smallmatrix}
            \right]~}A'\oplus C\xrightarrow{~\left[
              \begin{smallmatrix}
            y'x'&-c
              \end{smallmatrix}
            \right]~}C'\overset{\theta_P}{\dashrightarrow}$$
shows that the outside rectangle is a homotopy cartesian.
Moreover, the diagram (\ref{t2}) also shows that $\theta_R=x_{\ast}\theta_P$ holds.

Note that  the following diagram
$$
\xymatrix@C=1.6cm@R=1.4cm{A\ar@{=}[d]\ar[r]^{\left[
              \begin{smallmatrix}
          a\\ x
              \end{smallmatrix}
            \right]\qquad} &A'\oplus B\ar[r]^{\qquad\left[
              \begin{smallmatrix}
         x'&-b
              \end{smallmatrix}
            \right]\quad}\ar[d]^{\left[
              \begin{smallmatrix}
          1&0\\
          0&y
              \end{smallmatrix}
            \right]} &
B'\ar[d]^{y'}\ar@{-->}[r]^{\quad \theta_L}&\\
 A\ar[r]^{\left[
              \begin{smallmatrix}
        a\\yx
              \end{smallmatrix}
            \right]\qquad} & A'\oplus  C\ar[r]^{\quad\left[
              \begin{smallmatrix}
            y'x'&-c
              \end{smallmatrix}
            \right]} & C'\ar@{-->}[r]^{\quad \theta_P}&}
$$
is a morphism of $\E$-triangles,
this shows that $\theta_L=(y')^{\ast}\theta_P$ holds.

This completes the proof.  \qed

\begin{remark}
Our main result was
proved in \cite{M}, when $\C$ is an abelian category.
It also was proved in \cite{CF}, when $\C$ is an triangulated category.
\end{remark}

\textbf{Jing He}\\
College of Science, Hunan University of Technology and Business, 410205 Changsha, P. R. China\\
E-mail: jinghe1003@163.com
\\[0.3cm]
\textbf{Chenbei Xie}\\
College of Mathematics, Hunan Institute of Science and Technology, 414006 Yueyang, Hunan, P. R. China.\\
E-mail: xcb19980313@163.com
\\[0.3cm]
\textbf{Panyue Zhou}\\
College of Mathematics, Hunan Institute of Science and Technology, 414006 Yueyang, Hunan, P. R. China.\\
E-mail: panyuezhou@163.com

\end{document}